\documentclass[12pt,leqno]{amsart}
\usepackage{amssymb}
\usepackage{amsmath,amssymb,color}
\usepackage[none]{hyphenat}
\usepackage{geometry}
\numberwithin{equation}{section}

\newcommand{\la}{\lambda}
\newcommand{\va}{\varphi}
\newcommand{\ppp}{\partial}

\newcommand{\www}{\widetilde}

\newcommand{\R}{\mathbb{R}}
\newcommand{\C}{\mathbb{C}}
\newcommand{\N}{\mathbb{N}}

\newcommand{\sumn}{\sum_{n=1}^{\infty}}

\newcommand{\ooo}{\overline}
\newcommand{\OOO}{\Omega}

\newcommand{\MLO}{E_{\alpha,1}}
\newcommand{\MLT}{E_{\alpha,2}}
\newcommand{\UUU}{u(\alpha, p, h, H, \www{a})}
\newcommand{\VVV}{u(\beta, q, j, J, \www{b})}

\allowdisplaybreaks

\title[Simultaneous uniqueness for multiple parameters]
{Simultaneous uniqueness for multiple parameters identification in a fractional diffusion-wave equation}

\author{$^1\;$Xiaohua Jing, $^2\;$$^3\;$$^4\;$Masahiro Yamamoto}

\thanks{
$^1\;$School of Mathematics and Statistics,
Xi'an Jiaotong University,
Xi'an, 710049, China.
E-mail: {\tt xhjing5@163.com}
\\
$^2\;$Graduate School of Mathematical Sciences, The University
of Tokyo, Komaba, Meguro, Tokyo 153-8914, Japan, \\
$^3\;$Honorary Member of Academy of Romanian Scientists,
Splaiul Independentei Street, no 54,
050094 Bucharest Romania, \\
$^4\;$Peoples' Friendship University of Russia
(RUDN University) 6 Miklukho-Maklaya St, Moscow, 117198, Russian Federation.
E-mail: {\tt myama@ms.u-tokyo.ac.jp}
}

\date{}

\begin{document}

\maketitle

\baselineskip 18pt

\begin{abstract}
\sloppy{}
This article deals with the uniqueness in identifying multiple parameters
simultaneously in the one-dimensional time-fractional diffusion-wave equation
of fractional time-derivative order $\in (0,2)$
with the zero Robin boundary condition.
Using the Laplace transform and a transformation formula,
we prove the uniqueness in determining an order of the fractional derivative,
a spatially varying potential, initial values and Robin coefficients
simultaneously by boundary measurement data, provided that all the eigenmodes
of an initial value do not vanish.
Furthermore, for another formulation of inverse problem with input
source term in place of initial value, by the uniqueness in the case of
non-zero initial value and a Duhamel principle, we prove
the simultaneous uniqueness in determining multiple parameters for
a time-fractional diffusion-wave equation.
\end{abstract}



\maketitle

\section{Introduction and the main results}
\sloppy{}

In this article, we consider an inverse problem for the following one-dimensional time-fractional diffusion-wave equation firstly:
$$
\partial_{t}^{\alpha}u(x,t) = \partial_{x}^{2}u(x,t) - p(x)u(x,t), \quad \quad 0<x<1, \, 0<t<T, \eqno{(1.1)}
$$
$$
\partial_{x}u(0,t)-hu(0,t) =0, \  \partial_{x}u(1,t)+Hu(1,t)=0,  \quad  0<t<T, \eqno{(1.2)}
$$
$$
u(x,0) = a(x), \quad 0<x<1,  \eqno{(1.3)}
$$
$$
\partial_{t}u(x,0) = a^0(x), \quad 0<x<1, \quad if \ 1<\alpha<2.  \eqno{(1.4)}
$$
Here we assume that $\alpha \in (0,2)$, $p \in C^1[0,1]$,
$p\ge 0$ on $[0, 1]$ and $h, H>0$.
We use the Caputo fractional derivative of order
$\alpha$ for variable  $t$ defined by
\begin{eqnarray*}
\partial_{t}^{\alpha}u(t)=\frac{1}{\Gamma(n-\alpha)}\int ^{t}_{0}(t-s)^{n-\alpha-1}\frac{d^n u}{ds^n}(s)ds,
\end{eqnarray*}
if $\alpha >0$ satisfies $n-1<\alpha<n$ with $n\in \N$.
Here $\Gamma(\cdot)$ stands for the Gamma function.
The fractional order $\alpha$ is related to
the parameter specifying the large-time behavior of the waiting-time
distribution function \cite{Nig1986,Pod1998}.
Note that if the fractional order $\alpha = 1, 2$, then the equation
represents parabolic and hyperbolic equations, respectively.
Fractional diffusion-wave equations have been proposed for describing
for example,
anomalous diffusion phenomena characterized by the long-tailed profile in the
spatial distribution of densities as time passes.
Equation (1.1)-(1.4) is a model equation for anomalous diffusion in
heterogeneous media
and for related calculus and physical backgrounds, we can refer to
\cite{Met2000,Pod1998,Zhou2003}.

Most works on the research of the forward problem and inverse problem for the time-fractional diffusion-wave equations, we refer
to \cite{BW2012,BangtiIP2015, LYZZ, Yuri2013, SY, TWeiIP2016,Zhang2011}.
This research field is growing rapidly and we do not here give a complete
list of works.  In particular, as for inverse problems of determining
multiple parameters in the time-fractional diffusion equation simultaneously,
we can refer to \cite{JinIP2009,KYY2018,Li2013,LiIP2015,SunANM2017}.
However, most published works consider two kinds of parameters,
such as potential and fractional order.
On the other hand, for identifying multiple parameters including orders, to the authors' best knowledge,
the published paper focus on $0<\alpha<1$.
For $1<\alpha < 2$, there are very few publications, and
we can refer to \cite{LW2019,M2013IP}.
In the case where unknown orders $\alpha$ vary over $(0,2)$ including $1$,
there seem no theoretical results on the determination of orders and other
parameters.
The parameters $\alpha, p(x), a(x), a^0(x), h, H$ in (1.1)-(1.4) for the case $\alpha \in (0, 2)$ characterize physical properties of the diffusion process.
From the physical viewpoint, in modelling, it is not
natural that we a priori assume $0<\alpha\le 1$ and $1<\alpha<2$ separately,
and so it is a more feasible formulation of the inverse problem of determining
an order within $\alpha \in (0,2)$, restricted to neither $\alpha \in (0,1]$
nor $\alpha \in (1,2)$.

In this article,
we consider the inverse problem of uniqueness in the simultaneous identification of the fractional order derivative, potential, initial value and Robin coefficients in the boundary condition simultaneously for the model (1.1)-(1.4) from two boundary measurement data.

For the statement of our main results, we need to introduce some general
settings and notations firstly.  We assume that all the functions are
real-valued. We write
$$
(f,g) = \int^1_0 f(x)g(x) dx, \quad
\Vert f\Vert = (f,f)^{\frac{1}{2}}, \quad f, g \in L^2(0,1).
$$
Let $L^2(\Omega)$ be
a usual $L^2$-space with the inner product $( \cdot , \cdot)$ and $H^2$
denote the usual Sobolev spaces.

Then we define an operator $A_{p,h,H}$ by
$$\left\{ \begin{array}{rl}
&(A_{p,h,H} u)(x) = - u''(x) + p(x)u(x), \quad 0<x<1,\\
& \mathcal{D} (A_{p,h,H}) = \left\{
u \in H^2(0,1);\,  \frac{du}{dx}(0) - hu(0) = \frac{du}{dx}(1)
+ Hu(1) = 0 \right\}.
\end{array}\right.       \eqno{(1.5)}
$$
Let $\{ \la_n , \va_n \}_{n\in \N}$ and $\{\mu_n, \psi_n \}_{n\in \N}$
be the eigenvalues and
the eigenfunctions of the operators $A_{p,h,H}$ and
$A_{q,j,J}$ on the interval $[0,1]$, respectively:
$$\left\{ \begin{array}{rl}
& - \frac{d^2}{dx^2}\va_n(x) + p(x)\va_n(x) = \la_n\va_n(x), \quad
0<x<1, \\
& \frac{d}{dx}\va_n(0) - h\va_n(0) = \frac{d}{dx}\va_n(1) + H\va_n(1) = 0,\\
& \va_n(0) = 1,
\end{array}\right.
$$
and
$$\left\{ \begin{array}{rl}
& - \frac{d^2}{dx^2}\psi_n(x) + q(x)\psi_n(x) = \mu_n\psi_n(x), \quad
0<x<1, \\
& \frac{d}{dx}\psi_n(0) - j\psi_n(0) = \frac{d}{dx}\psi_n(1)
+ J\psi_n(1) = 0,\\
& \psi_n(0) = 1.
\end{array}\right.
$$
We set
$$
\rho_n = \Vert \va_n\Vert^2, \quad \sigma_n = \Vert \psi_n\Vert^2, \quad
n\in \N.
$$
Moreover, we choose
$$
(\va_n, \va_m) = \rho_n\delta_{nm} :=
\left\{ \begin{array}{rl}
\rho_n, \quad & n=m,\\
0, \quad & n\ne m
\end{array}\right.
$$
and $(\psi_n, \psi_m) = \sigma_n\delta_{nm}$.  Then it is known that
each of $\{ \va_n \}_{n\in \N}$ and $\{\psi_n \}_{n\in \N}$
is an orthogonal basis in $L^2(0,1)$.

Since the formulation of initial conditions changes according to the orders
$(0,1]$ and $(1,2)$, for convenience we introduce the following notations:
$$
\www{a} =
\left\{\begin{array}{rl}
a, \quad &\mbox{if $0<\alpha\le 1$},\\
(a, a^0), \quad &\mbox{if $1< \alpha < 2$},
\end{array}\right.
\quad
\www{b} =
\left\{\begin{array}{rl}
b, \quad &\mbox{if $0<\beta\le 1$},\\
(b, b^0), \quad &\mbox{if $1< \beta < 2$},
\end{array}\right.
$$
where $a, a^0 \in \mathcal{D}(A_{p,h,H})$ and
$b, b^0 \in \mathcal{D}(A_{q,j,J})$.

By $u(\alpha, p, h, H, \www{a})(x, t)$, we denote the solution to
(1.1)-(1.4),
in order to indicate its dependence on the parameters $\{\alpha, p, h, H,
\www{a}\}$.
Furthermore, for $n\in \N$, we define
$$
\www{a}^n :=
\left\{ \begin{array}{rl}
(a, \va_n), \quad & \mbox{if $0<\alpha\le 1$}, \\
\vert (a, \va_n)\vert + \vert (a^0, \va_n)\vert, \quad
&\mbox{if $1<\alpha < 2$}
\end{array}\right.
$$
and
$$
\www{b}^n :=
\left\{ \begin{array}{rl}
(b, \psi_n), \quad & \mbox{if $0<\beta\le 1$}, \\
\vert (b, \psi_n)\vert + \vert (b^0, \psi_n)\vert, \quad
&\mbox{if $1<\beta < 2$}.
\end{array}\right.
$$

Now we pose the main assumption in this article.

\

{\bf Assumption.}\\
{\it
Initial values $\www{a}$ and $\www{b}$ satisfy
$$
\vert \www{a}^n\vert + \vert \www{b}^n\vert \ne 0 \quad
\mbox{for all $n\in \N$}.
$$
}

Let $T>0$ be arbitrarily fixed.

Now we present our first main result on the uniqueness in simultaneously
determining multiple parameters in (1.1)-(1.4).
\\
\\
{\bf Theorem 1.} {\it
Let $0<\alpha, \beta<2$, $p, q \in C^1[0, 1]$, $p, q \ge0$
on $[0,1]$, $h, H, j, J >0$. Under {\bf Assumption}, if
$$
u(\alpha, p, h, H, \www{a})(\ell, t) = u(\beta, q, j, J, \www{b})
(\ell, t), \quad \ell=0,1, \,
\, 0<t<T,
$$
then
$$
\alpha=\beta, \ p(x)=q(x), \ \www{a}(x) = \www{b}(x), \quad 0<x<1,
\ \ and \ \ h=j, \, H=J.
$$
}

In the conclusion, we recall that $\www{a} = \www{b}$ means that
$a=b$ if $0<\alpha\le 1$ and $a=b, a^0=b^0$ if $1<\alpha<2$.

{\bf Assumption} requires the condition
$\vert \www{a}^n\vert + \vert \www{b}^n\vert \ne 0$ for all
$n\in \N$, which is related to unknown quantities $\alpha, \beta$ and
$\www{a}, \www{b}$.  We can make another interpretation:
Assuming that $(\alpha, p, h,H,\www{a})$ are known,
we are requested to identify $(\beta, q, j,J,\www{b})$ compared with
data from the known system with $(\alpha, p, h,H,\www{a})$.
Then we can replace the condition
$\vert \www{a}^n\vert + \vert \www{b}^n\vert \ne 0$ for all
$n\in \N$ by a condition on known quantities:
$$
\vert \www{a}^n\vert \ne 0 \quad \mbox{for all $n\in \N$}.
$$
Similar remarks hold for the second main result.
\\

{\bf Assumption} means that all the eigenmode of initial values
should be non-zero, and is a quite
restrictive condition.  As is seen by the proof, we can modify the
assumption, which allows us to choose an $N$-number inputs of
initial values whose union contains non-zero eigenmodes in all the eigenspaces.
More precisely,
\\
{\bf Theorem 1'}.
{\it
We assume that there exist $N\in \N$, and
$$
\www{a}_k:=
\left\{ \begin{array}{rl}
a_k\in \mathcal{D}(A_{p,h,H}), \, &\mbox{if $0<\alpha\le 1$},\\
(a_k, \, a_k^0) \in (\mathcal{D}(A_{p,h,H}))^2, \, &\mbox{if $1<\alpha < 2$},
\end{array}\right.
\quad
\www{b}_k:=
\left\{ \begin{array}{rl}
b_k\in \mathcal{D}(A_{q,j,J}), \, &\mbox{if $0<\beta\le 1$},\\
(b_k, \, b_k^0) \in (\mathcal{D}(A_{q,j,J}))^2, \, &\mbox{if $1<\beta < 2$},
\end{array}\right.
$$
for $k\in \{1, ..., N\}$ such that
$$
\bigcup_{k=1}^N \{ n\in \N;\, \vert \www{a}_k^n\vert
+ \vert \www{b}_k^n\vert \ne 0 \} = \N.
$$
Then
$$
u(\alpha,p,h,H,\www{a}_k)(\ell,t) = u(\beta,q,j,J,\www{b}_k)(\ell,t),
\quad \ell = 0,1, \, k \in \{1, ..., N\}, \, 0<t<T
$$
implies
$$
\left\{\begin{array}{rl}
&\alpha = \beta, \quad p(x) = q(x), \quad \www{a}_k(x) = \www{b}_k(x), \quad
0 \le x \le 1, \, k\in \{1, ..., N\},\\
& h=j, \quad H=J.
\end{array}\right.
$$
}
\\

Here we write
$$
\www{a}_k^n :=
\left\{ \begin{array}{rl}
(a_k, \va_n), \quad &\mbox{if $0<\alpha\le 1$},\\
\vert (a_k, \va_n)\vert + \vert (a_k^0, \va_n)\vert, \quad
&\mbox{if $1<\alpha < 2$},
\end{array}\right.
$$
and
$$
\www{a}_k^n:=
\left\{ \begin{array}{rl}
(b_k, \psi_n), \quad &\mbox{if $0<\beta\le 1$},\\
\vert (b_k, \psi_n)\vert + \vert (b_k^0, \psi_n)\vert, \quad
&\mbox{if $1<\beta < 2$},
\end{array}\right.
$$
for $k\in \{1, ..., N\}$ and $n\in \N$.

In this article, we are interested also in the following inverse problem
for a time-fractional diffusion equation with source term:
$$
\partial_{t}^{\alpha}\www{u}(x,t) = \partial_{x}^{2}\www{u}(x,t)- p(x)\www{u}(x,t) + \theta(t)g(x) ,  \quad 0<x<1, \, 0<t<T, \eqno{(1.6)}
$$
$$
\partial_{x}\www{u}(0,t)-h\www{u}(0,t)=0, \ \partial_{x}\www{u}(1,t)+H\www{u}(1,t)=0,  \ 0<t<T,\eqno{(1.7)}
$$
$$
\www{u}(x,0) =0, \quad 0<x<1,  \eqno{(1.8)}
$$
$$
\partial_{t}\www{u}(x,0) = 0, \ 0<x<1, \quad if \ 1<\alpha<2.  \eqno{(1.9)}
$$

Equations (1.6)-(1.9) are closely
related to equations (1.1)-(1.4) by using a fractional Duhamel principle \cite{GYM,LRY,LYZZ}.
Theorem 1 establishes the uniqueness
for the multiple parameters simultaneously in equation (1.1)-(1.4).
Therefore, taking advantage of the fractional Duhamel principle, we can
transfer the uniqueness by Theorem 1 into the uniqueness for
(1.6)-(1.9) by using the boundary measurement data.
Here and henceforth let $\widetilde{u}(p, h, H)$ denote
the solution to (1.6)-(1.9) with $\{p, h, H\}$ and $g\in H^2(0,1)$.
Moreover, we fix
$$
\theta \in C^1[0,T], \quad \not\equiv 0.
$$
Then\\
\\
{\bf Theorem 2.}
{\it
Let $0<\alpha, \beta<2$, $p, q \in C^1[0,1]$, $p, q \ge 0$ on $[0,1]$,
$h, H,  j, J >0$.
Under {\bf Assumption}, if
$$
\widetilde{u}(p, h, H)(\ell, t) = \widetilde{u}(q, j, J)(\ell, t), \quad \ell=0,1, \
\, 0<t<T,
$$
then
$$
p(x)=q(x), \quad 0<x<1, \
\ and \ \ h=j, \, H=J.
$$
}

The rest of the article is composed of three sections and an appendix.
In Section 2, we show preliminary results on the Mittag-Leffler function
and provide some lemmas that are used for the proofs of the main results.
In Section 3, we complete the proof of Theorem 1.
Based on Theorem 1 and the fractional Duhamel
principle, we complete the proof of Theorem 2 in Section 4.
The appendix is devoted to the proof of Lemma 4 in Section 2.
Throughout the article, we denote by $C$ a generic constant,
which may differ at different occurrences.

\section{Preliminaries}
\sloppy{}

In this section, we recall some preliminary results on the Mittag-Leffler function and
provide some lemmas which are needed in our subsequent arguments.

To start with, we recall the two-parameter Mittag-Leffler function
$E_{\alpha,\beta}(z)$:
$$
E_{\alpha,\beta}(z)=  \sum_{k=0}^{\infty}\frac{z^{k}}{\Gamma(k \alpha + \beta)} \quad z\in \C, \eqno{(2.1)}
$$
where $\alpha > 0$ and $\beta \in \mathbb{R}$ are arbitrary constants
(e.g., \cite{Pod1998}). By the power series, $E_{\alpha,\beta}(z)$ is an entire function of $z\in \C$. Note that $E_{1,1}(z)=e^z.$
Moreover, we state the following property on the Mittag-Leffler function.\\
\\
{\bf Lemma 1 (\cite{Pod1998}).}
{\it
Let $0<\alpha<2$ and $\beta>0$ be arbitrary.
We suppose that $\mu$ satisfies $\pi \alpha/2<\mu< \min \{\pi, \pi \alpha\}$.
Then there exists a constant $C = C(\alpha, \beta, \mu)>0$
such that
$$
|E_{\alpha,\beta}(z)|\leq \frac{C}{1+|z|}, \quad \quad z\in\C, \quad \mu \leq
|\mbox{arg} \ z|\leq \pi.
$$
}

From the definition of operator $A_{p,h,H}$ in (1.5), for all $n \geq 0$
it is known that
$$
\sqrt{\lambda_{n}} = n\pi + \frac{\omega}{n} + \frac{k_{n}}{n}, \ \ \ \ \
\sumn \vert k_n\vert^2 < \infty,     \eqno{(2.2)}
$$

$$
\varphi_{n}(x;\lambda_{n}) =\sqrt{2} cos(n \pi x) + \frac{\xi_{n}(x)}{n}, \ \ \ \ \  |\xi_{n}(x)| \leq C, \eqno{(2.3)}
$$
where $C>0$ is a constant and
$\omega = h + H + \frac{1}{2}\int_{0}^{1}q(t)dt$
(e.g., Levitan and Sargsjan \cite{LS}).

\

If $a \in \mathcal{D}(A_{p,h,H})$, then $\sumn \vert \frac{(a,\va_n)}{\rho_n} \vert < \infty$.
Indeed
$$
(a, \va_n) = \frac{1}{\la_n}(a, \la_n\va_n)
= \frac{1}{\la_n}(a, A_{p,h,H}\va_n) = \frac{1}{\la_n}(A_{p,h,H}a, \va_n).
$$
Hence
$$
\vert (a, \va_n) \vert \le \frac{1}{\la_n}\Vert A_{p,h,H}a\Vert
\Vert \va_n\Vert \le \frac{C}{n^2}\Vert A_{p,h,H}a\Vert \sqrt{\rho_n} ,
$$
which implies
$$
\sumn \left\vert \frac{(a,\va_n)}{\rho_n} \right\vert
\le  C\sumn \frac{1}{n^2}\frac{1}{\sqrt{\rho_n}} \Vert A_{p,h,H}a\Vert
< \infty.
$$
Henceforth we set
$$
p_n = \frac{(a,\va_n)}{\rho_n}, \ p^0_n = \frac{(a^0,\va_n)}{\rho_n},
\ q_n = \frac{(b,\psi_n)}{\sigma_n}, \ q^0_n = \frac{(b^0,\psi_n)}{\sigma_n},
\quad n\in \N.
$$
Then, we can see
$$
\sumn (\vert p_n\vert + \vert p_n^0\vert +\vert q_n\vert
+ \vert q_n^0\vert) < \infty.        \eqno{(2.4)}
$$

\

Furthermore, we show some lemmas.
In particular, the first two lemmas are concerned with a
formula, connecting eigenfunctions of the spatial difference operator through an integral transformation.
The proof can be found in \cite{TS1985} for example.
Let $\OOO = \{ (x,y);\, 0<y<x<1\}$.

\

{\bf Lemma 2.}
{\it For each given $p, q \in C^1[0,1]$ and $h, j\in \R$, there exists a unique $K = K(x,y)= K(x,y;p,h; q,j) \in C^2(\ooo{\OOO})$ such that
$$
\left\{ \begin{array}{rl}
& \ppp_x^2K(x,y) - \ppp_y^2K(x,y) = q(x)K(x,y) - p(y)K(x,y) \quad
(x,y) \in \ooo{\OOO}, \\
& \ppp_yK(x,0) = hK(x,0) \quad 0\leq x \leq 1,
\end{array}\right.
$$
and
$$
K(x,x) = j-h + \frac{1}{2}\int^x_0 (p(\xi) - q(\xi)) d\xi, \quad
0 \leq x \leq 1.                                        \eqno{(2.5)}
$$
}

\

{\bf Lemma 3 (Transformation formula).}
{\it
Let $K$ be defined in Lemma 2, and let $\va(x,\la)$ satisfy
$$
\left\{ \begin{array}{rl}
& -\frac{d^2\va}{dx^2}(x) + p(x)\va(x) = \la \va(x), \quad 0 \leq x \leq 1, \\
& \frac{d\va}{dx}(0) = h, \quad \va(0) = 1
\end{array}\right.
$$
Then
$$
\psi(x,\la) = \va(x,\la)
+ \int^x_0 K(x,y)\va(y,\la) dy, \quad 0 \leq x \leq 1, \quad \la \in \R,
$$
satisfies
$$
\left\{ \begin{array}{rl}
& -\frac{d^2\psi}{dx^2}(x) + q(x)\psi(x) = \la \psi(x), \quad 0 \leq x \leq 1, \\
& \frac{d\psi}{dx}(0) = j, \quad \psi(0) = 1.
\end{array}\right.
$$
}

We conclude this section with the following lemma whose proof is
given in Appendix.
\\
\\
{\bf Lemma 4.} {\it
Let $\UUU$ satisfy (1.1)-(1,4), then $\UUU(0,t) = 0$ for $0<t<T$ implies $\UUU(x,t) = 0$ for
$0<x<1$ and $0<t<T$.
}

\section{Proof of Theorem 1}
\sloppy{}

The proof is divided into Step I and Step II.
In Step I, we prove the uniqueness of fractional order and
establish the uniqueness of other quantities in Step II.

Note that $E_{1,1}(z)=e^z.$
We have
$$
\left\{\begin{array}{rl}
& \UUU(x,t) = \sum_{n=1}^\infty p_n E_{\alpha,1}(-\la_nt^{\alpha})\va_n(x),
\cr\\
& \VVV(x,t) = \sum_{n=1}^\infty q_n
E_{\beta,1}(-\mu_nt^{\beta})\psi_n(x)
\end{array}\right.                \eqno{(3.1)}
$$
for $\alpha, \beta \in (0,1]$, and
$$
\left\{\begin{array}{rl}
&\UUU(x,t) = \sum_{n=1}^\infty p_n
E_{\alpha,1}(-\la_nt^{\alpha})\va_n(x) + \sum_{n=1}^\infty p^0_n
tE_{\alpha,2}(-\la_nt^{\alpha})\va_n(x),\cr\\
&\VVV(x,t) = \sum_{n=1}^\infty q_n
E_{\beta,1}(-\mu_nt^{\beta})\psi_n(x)
+ \sum_{n=1}^\infty q^0_n tE_{\beta,2}(-\mu_nt^{\beta})\psi_n(x)
\end{array}\right.               \eqno{(3.2)}
$$
for $\alpha, \beta \in (1,2)$,
which are proved in e.g., \cite{SY} and $\alpha, \beta\ne 1$ and the case
$\alpha=1$ and $\beta=1$ is classical.
We recall
$$
p_n = \frac{(a,\va_n)}{\rho_n}, \ p^0_n = \frac{(a^0,\va_n)}{\rho_n},
\ q_n = \frac{(b,\psi_n)}{\sigma_n}, \ q^0_n = \frac{(b^0,\psi_n)}{\sigma_n},
\quad n\in \N.
$$
By {\bf Assumption}, Lemma 1 and (2.2)-(2.4),
we see that
$\UUU(\ell,t)$ and $\VVV(\ell,t)$ with $\ell=0,1$ are convergent in
$C[\delta,\infty)$ with arbitrary $\delta>0$.
\\

Now, we are ready to prove Theorem 1.

\

{\bf Step I.} Firstly, we will prove $\alpha=\beta$.
We can exclude $\alpha=\beta=1$ trivially.

We apply an argument in Liu, Hu and Yamamoto \cite{LHY2020}, which
relies on the analysis of the poles of Laplace transformed data.

We can classify all the cases into the following three cases:

\

\begin{itemize}
 \item {\bf Case 1.}
$(\alpha, \beta) \in (0,1]^2\setminus \{(1,1)\}$.
 \item {\bf Case 2.}
$0<\alpha\leq1$ and $1<\beta<2$ or $1<\alpha<2$ and $0<\beta\leq1$.
\item
{\bf Case 3.} $1<\alpha, \beta<2$.
\end{itemize}

\

{\bf Case 1:} $(\alpha, \beta) \in (0,1]^2\setminus \{(1,1)\}$.

Without loss of generality, we can assume that $\alpha>\beta$.
Since
$$
\UUU(0,t) = \VVV(0,t), \quad 0<t<T,
$$
by Lemma 1 and (2.2)-(2.4), the analyticity in $t>0$ (\cite{SY}) yields
$$
\sumn p_{n}E_{\alpha,1}(-\la_nt^{\alpha})
= \sumn q_{n} E_{\beta,1}(-\mu_nt^{\beta}), \quad t > 0. \eqno{(3.3)}
$$
Taking the Laplace transform in terms of (3.3) and using the Laplace transform
of $E_{\alpha,1}$ (e.g., \cite{Pod1998}), we see
$$
\sumn \frac{p_{n}z^{\alpha-1}}{z^{\alpha}+\la_n}
= \sumn \frac{q_{n}z^{\beta-1}}{z^{\beta}+\mu_n}, \quad \mbox{Re}\, z > C_0,
$$
where $C_0>0$ is a constant.

We choose large $\ell>0$ such that $\alpha\ell$, $\beta\ell>1$. Setting
$z = \eta^\ell$, $\gamma=\alpha \ell$ and $\delta=\beta \ell$,
we obtain
$$
\sumn \frac{p_n \eta^{\gamma}}{\eta^\gamma +\la_n}
=
\sumn \frac{q_n \eta^{\delta}}{\eta^\delta + \mu_n}, \ \ \eta \in Q.
$$
We set
$$
Q := \{ \eta\in \C; \eta\neq0, | \mbox{arg}\, \ \eta|< \frac{1}{2}\pi \}
\cap \{\C \backslash(\{ \eta^\pm_{m,\gamma}\}_{m\in \N} \cup \{ \eta^\pm_{n,\delta}\}_{n\in \N})\},
$$
and
$$
\eta^\pm_{m,\gamma} := \la^{\frac{1}{\gamma}}_m e^{ \pm \sqrt{-1} \frac{\pi}{\gamma}}, \quad
\eta^\pm_{n,\delta} := \mu^{\frac{1}{\delta}}_n e^{\pm \sqrt{-1} \frac{\pi}{\delta}}.
$$
Since $\alpha>\beta$, we have $\gamma>\delta$. Hence
$$
\sumn \frac{p_n \eta^{\gamma-\delta}}{\eta^\gamma +\la_n}
= \sumn \frac{q_n}{\eta^\delta + \mu_n},   \ \ \ \ {\eta\in Q}
$$
For arbitrary $n_0 \in \N$, we obtain
$$
\frac{p_{n_0} \eta^{\gamma-\delta}}{\eta^\gamma +\la_{n_0}}
= - \sum_{n\neq n_0} \frac{p_n \eta^{\gamma-\delta}}{\eta^\gamma +\la_n}
+ \sumn \frac{q_n }{\eta^\delta + \mu_n}.
$$
We rewrite $J(\eta)=- \sum_{n\neq n_0} \frac{p_n \eta^{\gamma-\delta}}{\eta^\gamma +\la_n}
+ \sumn \frac{q_n}{\eta^\delta + \mu_n}$.
Thus,
$$
p_{n_0} \eta^{\gamma-\delta} = (\eta^\gamma + \la_{n_0})J(\eta), \quad
\eta\in Q.                           \eqno{(3.4)}
$$
Set
$$
\eta_{n_0} := \la^{\frac{1}{\gamma}}_{n_0} e^{ \sqrt{-1} \frac{\pi}{\gamma}}.
$$
Then,
$$
\lim_{\eta\in Q,\eta\to \eta_{n_0}} J(\eta)
$$
exists.
Indeed, for $ n \neq n_0$,
$$
(\eta_{n_0})^\gamma = \la_{n_0}e^{ \sqrt{-1} \pi} = - \la_{n_0} \neq -\la_n.
$$
Moreover, if $\la_{n_0}^{ \frac{\delta}{\gamma}} e^{ \sqrt{-1}
\pi \frac{\delta}{\gamma}} = -\mu_{m_0}$
with some $m_0 \in \N$, then $\la_n >0$ and $\mu_n >0$ yield
$e^{ \sqrt{-1} \pi \frac{\delta}{\gamma}}
= e^{ \sqrt{-1} \pi \frac{\beta}{\alpha}} < 0$.  Therefore
$\frac{\beta}{\alpha} = \pi,\ \ 3\pi, \cdots$, which is impossible
by $\alpha>\beta.$  Thus
$$
(\eta_{n_0})^\delta = \la_{n_0}^{ \frac{\delta}{\gamma}}  e^{ \sqrt{-1} \pi \frac{\delta}{\gamma}} \not\in \{-\mu_n \}_{n\in \N}.
$$
Hence we see that $\lim_{\eta\in Q,\eta\to \eta_{n_0}} J(\eta)
= J(\eta_{n_0})$.

Letting $ \eta \rightarrow \eta_{n_0}$ in (3.4), we have
$p_{n_0}(\eta_{n_0})^{\gamma-\delta}=0$,
that is, $p_{n_0}=0$.  Since $n_0$ ia arbitrary,
we have that $(a, \va_n)=0$ for all $n\in \N$.  Hence $\UUU(x,t) = 0$ for
$0<x<1$ and $0<t<T$.  By $\UUU(0,t) = \VVV(0,t)$ for $0<t<T$, Lemma 4
yields $\VVV(x,t) = 0$ for $0<x<1$ and $0<t<T$.
Hence $q_n=0$ for all $n\in \N$.
Thus $p_n=q_n=0$ for all $n\in \N$, which implies
$\vert \www{a}^n\vert + \vert \www{b}^n\vert = 0$ for all
$n\in \N$.  This is a contradiction to
{\bf Assumption}.
Hence, we see that $\alpha > \beta$ is impossible. Similarly, $\beta>\alpha$ is impossible. Therefore $\alpha=\beta$ follows.

\

{\bf Case 2:} $0<\alpha\leq1$ and $1<\beta<2$, or $1<\alpha<2$ and $0<\beta\leq1$.

Here, we prove the case of $0<\alpha\leq1$ and $1<\beta<2$, and the proof in the case of $1<\alpha<2$ and $0<\beta\leq1$ is similar.
Since we have
$$
\UUU(0,t) = \VVV(0,t), \quad 0<t<T,
$$
then the $t$-analyticity (\cite{SY}) yields
$$
\sumn p_{n}E_{\alpha,1}(-\la_nt^{\alpha})
= \sumn q_{n} E_{\beta,1}(-\mu_nt^{\beta})
+  \sumn q^0_{n}t E_{\beta,2}(-\mu_nt^{\beta})    , \quad t > 0. \eqno{(3.5)}
$$
We choose large $\ell>0$ such that $\alpha\ell>1$, $\beta\ell>2$, and $(\beta-\alpha)\ell>1$. Setting
$z = \eta^\ell$, $\gamma=\alpha \ell$, $\delta=\beta \ell$.
By the Laplace transform in terms of (3.5), we have
$$
\sumn \frac{p_n \eta^{\gamma}}{\eta^\gamma +\la_n}
=
\sumn \frac{q_n \eta^{\delta} + q^0_n \eta^{\delta-1}}{\eta^\delta + \mu_n}, \ \ \eta \in Q.
$$
Set
$$
Q := \{ \eta \in \C; \eta \neq 0, \ |arg \ \eta|< \pi-\varepsilon \}
\cap \{\C \backslash(\{ \eta^\pm_{m,\gamma}\}_{m\in \N} \cup \{ \eta^\pm_{n,\delta}\}_{n\in \N})\},
$$
where
$$
\eta^\pm_{m,\gamma} := \la^{\frac{1}{\gamma}}_m e^{ \pm \sqrt{-1} \frac{\pi}{\gamma}}, \quad
\eta^\pm_{n,\delta} := \mu^{\frac{1}{\delta}}_n e^{\pm \sqrt{-1} \frac{\pi}{\delta}}.
$$
Therefore, we have
$$
\sumn \frac{p_n}{\eta^\gamma +\la_n}
=
\sumn \frac{q_n \eta^{\delta-\gamma} + q^0_n \eta^{\delta-\gamma-1}}{\eta^\delta + \mu_n}.  \ \ \ \ {\eta\in Q}
$$
Let $n_0$ be arbitrary.  Then
$$
\frac{q_{n_0} \eta^{\delta-\gamma}+q^0_{n_0} \eta^{\delta-\gamma-1}}{\eta^\delta + \mu_{n_0}}
= -\sum_{n\neq n_0}\frac{q_n \eta^{\delta-\gamma}+q^0_n \eta^{\delta-\gamma-1}}{\eta^\delta + \mu_n}
+ \sumn \frac{p_n}{\eta^\gamma +\la_n}.
$$
We rewrite $J(\eta)=-\sum_{n\neq n_0}\frac{q_n \eta^{\delta-\gamma}+q^0_n \eta^{\delta-\gamma-1}}{\eta^\delta + \mu_n}
+ \sumn \frac{p_n}{\eta^\gamma +\la_n}.$
Thus,
$$
q_{n_0} \eta^{\delta-\gamma}+q^0_{n_0} \eta^{\delta-\gamma-1}
 = (\eta^\delta + \mu_{n_0})J(\eta), \quad \eta\in Q.   \eqno{(3.6)}
$$
Set
$$
\eta_{n_0} := \mu^{\frac{1}{\delta}}_{n_0} e^{ \sqrt{-1} \frac{\pi}{\delta}}.
$$
Then
$$
\lim_{\eta\in Q,\eta\to \eta_{n_0}} J(\eta)
$$
exists.
Indeed, for $ n \neq n_0$,
$$
(\eta_{n_0})^\delta = \mu_{n_0}e^{ \sqrt{-1} \pi} = - \mu_{n_0} \neq -\mu_n.
$$
Moreover,
if $\mu_{n_0}^{ \frac{\gamma}{\delta}} e^{ \sqrt{-1} \pi \frac{\gamma}{\delta}} = -\la_{m_0}$
with some $m_0 \in \N$, then $\la_n >0, \ \mu_n >0$ yield
$
e^{ \sqrt{-1} \pi \frac{\gamma}{\delta}} < 0.
$
Therefore $\frac{\gamma}{\delta}\pi = \frac{\alpha}{\beta}\pi = \pi,\ \ 3\pi, \cdots$, which is impossible by $\alpha\in(0,1]$,
$\beta\in(1,2).$ Thus
$$
(\eta_{n_0})^\gamma = \mu_{n_0}^{ \frac{\gamma}{\delta}}  e^{ \sqrt{-1} \pi \frac{\gamma}{\delta}} \not\in \{-\la_n \}_{n\in \N}.
$$
Hence we see that $\lim_{\eta\in Q,\eta\to \eta_{n_0}} J(\eta)
= J(\eta_{n_0})$.

Letting $ \eta \rightarrow \eta_{n_0}$ in (3.6), we have
$\eta_{n_0}^{\delta-\gamma-1}(q_{n_0}\eta_{n_0} + q^0_{n_0})=0$.
That is,
$q_{n_0}= q^0_{n_0}=0.$
Since $n_0$ ia arbitrary, we have that $(b, \psi_n)=(b^0, \psi_n)=0, n\in \N$.
Hence $\VVV(x,t) = 0$ for
$0<x<1$ and $0<t<T$.  By $\UUU(0,t) = \VVV(0,t)$ for $0<t<T$, Lemma 4
yields $\UUU(x,t) = 0$ for $0<x<1$ and $0<t<T$.
Hence $p_n=0$ for all $n\in \N$.
Thus $p_n=q_n=q^0_{n_0}=0$ for all $n\in \N$, which implies
$\vert \www{a}^n\vert + \vert \www{b}^n\vert = 0$ for all
$n\in \N$.  This is a contradiction to
{\bf Assumption}.
Therefore, we see that this case is impossible.

\

{\bf Case 3:} $1<\alpha, \beta<2$.

We have
$$
\UUU(0,t) = \VVV(0,t), \quad 0<t<T,
$$
the $t$-analyticity yields
$$
\sumn p_{n}E_{\alpha,1}(-\la_nt^{\alpha}) + \sumn p^0_{n}tE_{\alpha,2}(-\la_nt^{\alpha})
= \sumn q_{n} E_{\beta,1}(-\mu_nt^{\beta}) + \sumn q^0_{n}tE_{\beta,2}(-\mu_nt^{\beta}), \quad t > 0.
$$
By the Laplace transform, we have
$$
\sumn \frac{p_n z^{\alpha-1}+ p^0_n z^{\alpha-2}}{z^\alpha +\la_n}
=
\sumn \frac{q_n z^{\beta-1}+ q^0_n z^{\beta-2}}{z^\beta + \mu_n}, \ \ z\in Q.     \eqno{(3.7)}
$$
We set
$$
Q := \{ z\in \C; z\neq 0, \ |arg \ z|< \frac{\pi}{2} \}
\cap \{\C \backslash(\{ z^\pm_{m,\alpha}\}_{m\in \N} \cup \{ z^\pm_{n,\beta}\}_{n\in \N})\},
$$
where
$$
z^\pm_{m,\alpha} := \la^{\frac{1}{\alpha}}_n e^{\pm \sqrt{-1} \frac{\pi}{\alpha}}, \quad
z^\pm_{n,\beta} := \mu^{\frac{1}{\beta}}_n e^{\pm \sqrt{-1} \frac{\pi}{\beta}}.
$$
Without loss of generality, we assume that $\alpha>\beta$. Hence
$$
\sumn \frac{p_n z^{\alpha-\beta+1}+ p^0_n z^{\alpha-\beta}}{z^\alpha +\la_n}
=
\sumn \frac{q_n z+ q^0_n }{z^\beta + \mu_n}.  \ \ \ \ {z\in Q}
$$
Let $n_0$ be arbitrary. Then
$$
\frac{p_{n_0} z^{\alpha-\beta+1}+ p^0_{n_0} z^{\alpha-\beta}}{z^\alpha +\la_{n_0}}
= - \sum_{n\neq n_0} \frac{p_n z^{\alpha-\beta+1}+ p^0_n z^{\alpha-\beta}}{z^\alpha +\la_n}
+ \sumn \frac{q_n z+ q^0_n }{z^\beta + \mu_n}.
$$
We rewrite $J(z)=- \sum_{n\neq n_0} \frac{p_n z^{\alpha-\beta+1}+ p^0_n z^{\alpha-\beta}}{z^\alpha +\la_n}
+ \sumn \frac{q_n z+ q^0_n }{z^\beta + \mu_n}$.
Thus,
$$
p_{n_0} z^{\alpha-\beta+1}+ p^0_{n_0} z^{\alpha-\beta}= (z^\alpha + \la_{n_0})J(z), \quad z\in Q.   \eqno{(3.8)}
$$
Set
$$
z_{n_0} := \la^{\frac{1}{\alpha}}_{n_0} e^{ \sqrt{-1} \frac{\pi}{\alpha}}.
$$
Then,
$$
\lim_{z\in Q, z\to z_{n_0}} J(z)
$$
exists. Indeed, for $ n \neq n_0$, we know
$$
(z_{n_0})^\alpha = \la_{n_0}e^{\sqrt{-1} \pi} = - \la_{n_0} \neq -\la_n.
$$
Moreover,
if $\la_{n_0}^{ \frac{\beta}{\alpha}} e^{\sqrt{-1} \pi \frac{\beta}{\alpha}} = -\mu_{m_0}$
with some $m_0 \in \N$, then $\la_n >0, \ \mu_n >0$ yields
$ e^{\sqrt{-1} \pi \frac{\beta}{\alpha}} < 0 $.
Therefore $\frac{\beta}{\alpha} = \pi,\ \ 3\pi, \cdots$, which is impossible by $\alpha>\beta.$
Thus
$$
(z_{n_0})^\beta = \la_{n_0}^{ \frac{\beta}{\alpha}}  e^{ \sqrt{-1} \pi \frac{\beta}{\alpha}} \not\in \{-\mu_n \}_{n\in \N}.
$$
Hence we see that $\lim_{z\in Q, z\to z_{n_0}} J(z)
= J(z_{n_0})$.

Letting $ z \rightarrow z_{n_0}$ in (3.8), we have
$(z_{n_0})^{\alpha-\beta}(p_{n_0}z_{n_0} + p^0_{n_0})=0$.
That is,
$p_{n_0}= p^0_{n_0}=0.$
Since $n_0$ ia arbitrary,
we have that $(a, \va_n)=(a^0, \va_n)=0$ for all $n\in \N$.  Hence $\UUU(x,t) = 0$ for
$0<x<1$ and $0<t<T$.  By $\UUU(0,t) = \VVV(0,t)$ for $0<t<T$, Lemma 4
yields $\VVV(x,t) = 0$ for $0<x<1$ and $0<t<T$.
Hence $q_n=q^0_n=0$ for all $n\in \N$.
Thus $p_n=p^0_n=q_n=q^0_n=0$ for all $n\in \N$, which implies
$\vert \www{a}^n\vert + \vert \www{b}^n\vert = 0$ for all
$n\in \N$.  This is a contradiction to
{\bf Assumption}.
Hence, we see that $\alpha > \beta$ is impossible. Similarly, $\beta>\alpha$ is impossible. Therefore $\alpha=\beta$ follows.

Therefore, from the above arguments,
we obtain $\alpha=\beta$, $\alpha$, $\beta \in(0,2)$.
\\
\\
{\bf Step II.} Next, we prove the uniqueness of parameters $p, h, H,
\widetilde{a}$.

\

When $\alpha, \beta = 1$, the inverse problem is concerned with a
parabolic equation and
the uniqueness result has been given in \cite{TS1986, TSRM1980}.
Therefore, we consider the uniqueness of $p, h, H, \widetilde{a}$
in Theorem 1 separately in the following two cases.
\begin{itemize}
 \item {\bf{ Case 1.} $\alpha \in (0,1)$};
 \item {\bf{ Case 2.} $\alpha \in (1,2)$}.
\end{itemize}

\

{\bf{Case 1:} $\alpha \in (0,1)$.}

Since $\alpha=\beta$ is already proved,
$\UUU(0,t) = \VVV(0,t), \ 0<t<T$, with $t$-analyticity yields
$$
\sumn p_{n}E_{\alpha,1}(-\la_nt^{\alpha})
= \sumn q_{n} E_{\alpha,1}(-\mu_nt^{\alpha}), \quad \quad t > 0.
$$
We take the Laplace
transforms termwise of both sides of the above equation to obtain
$$
\sumn \frac{p_n z^{\alpha-1}}{z^\alpha +\la_n}
=
\sumn \frac{q_n z^{\alpha-1}}{z^\alpha + \mu_n}, \quad \quad Re \ z>0.
$$
That is
$$
\sumn \frac{p_n}{\xi +\la_n}
=
\sumn \frac{q_n }{\xi + \mu_n}, \quad \quad  Re \ \xi>0.  \eqno{(3.9)}
$$

By the asymptotic expression of the eigenvalues and the
eigenfunctions in \cite{LS},
we can analytically continue both sides of (3.9) in $\xi$ and
the above series are convergent uniformly on any compact set in
$\C \setminus(\{-\la_n\}_{n\geq 1} \cup \{-\mu_n\}_{n\geq 1})$.

Now we will prove that $\la_n = \mu_n$ for $n\in \N$.
Assume that there exists $n_0\in \N$ such that
$\la_{n_0} \neq \mu_m$ for all $m\in \N$.
Then we can take a suitable disk which includes
$-\la_{n_0}$ and does not include $\{ -\la_n\}_{n \neq n_0} \cup \{ -\mu_m\}_{m\geq 1}$.
Integrating (3.9) on the boundary of the disk, we obtain
$2\pi i p_{n_0}=0$. This is impossible because of $p_{n_0} \neq 0$.
Indeed, if $p_{n_0}=0$, since $n_0$ is arbitrary, then $p_{n_0}= 0$ yields $(a, \varphi_n)=0$ for all $n\in\N$.
Hence $u(p, h, H, \widetilde{a})(x,t) = 0$ for
$0<x<1$ and $0<t<T$.  By $u(p, h, H, \widetilde{a})(0,t) = u(q, j, J,\widetilde{b})(0,t)$ for $0<t<T$, Lemma 4
yields $u(q, j, J,\widetilde{b} )(x,t) = 0$ for $0<x<1$ and $0<t<T$.
Hence $q_n=0$ for all $n\in \N$.
Thus $p_n=q_n=0$ for all $n\in \N$, which implies
$\vert \www{a}^n\vert + \vert \www{b}^n\vert = 0$ for all
$n\in \N$.  This is a contradiction to
{\bf Assumption}.
Hence, we have $p_{n_0} \neq 0$.
Therefore, for each $n\in \N$, there exists $m(n)\in\N$ such that
$$
\la_n = \mu_{m(n)}, \quad n\geq 1.
$$
By the asymptotics (2.2) of the eigenvalues, we have
$$
\sqrt{\la_n}=n\pi + O\left(\frac{1}{n}\right), \ \ \sqrt{\mu_n}=n\pi
+ O\left(\frac{1}{n}\right)
$$
as $n\rightarrow \infty$, and so
$$
n\pi + O\left(\frac{1}{n}\right)=m(n)\pi + O\left(\frac{1}{m(n)}\right)
$$
as $n\rightarrow \infty$.
Consequently $\lim_{n\rightarrow \infty}(m(n)-n)=0$ by $|m(n)-n|=0$ or
$\geq 1$.
We can find $\widetilde{N}\in\N$ such that $m(n)=n$ for all
$n\geq \widetilde{N}$.
This implies that there exists exactly $m(\widetilde{N})$ eigenvalues $\mu_n$
on $[0, \mu_{m(\widetilde{N})}]$
and exactly $\widetilde{N}$ eigenvalues $\la_n$ on $[0, \la_{\widetilde{N}}]=[0, \mu_{m(\widetilde{N})}]$.
Hence $m(n)=n$ for each $n\in\N$. Then we
we can obtain
$$
\la_n = \mu_n, \quad n\geq 1.
$$
\\

Moreover, from $u(p, h, H, \widetilde{a})(0,t) = u(q, j, J,\widetilde{b})(0,t)$, we have
$$
\sumn p_n\MLO(-\la_nt^{\alpha})
=\sumn q_n\MLO(-\la_nt^{\alpha})
\quad \quad t>0.
$$
We choose large $\ell>0$ such that $\alpha\ell>1$. Setting
$z = \eta^\ell$, $\gamma=\alpha \ell$.
By the Laplace transform, we have
$$
\sumn \frac{p_n}{\eta^\gamma +\la_n}
=
\sumn \frac{q_n }{\eta^\gamma + \la_n}, \ \ Re \ \eta >0.
$$
Set
$$
\eta_{m, \gamma}^{\pm}= \la^{\frac{1}{\gamma}}_m e^{\pm \sqrt{-1} \frac{\pi}{\gamma}}.
$$
Letting $\eta \to \eta_{m, \gamma}^{\pm}$, we have
$$
p_m =q_m. \quad m\in\N. \eqno{(3.10)}
$$

Next, by applying the analyticity in $t>0$,
equation $u(p, h, H, \widetilde{a})(1,t) = u(q, j, J,\widetilde{b})(1,t)$,
$0<t<T$ yields
$$
\sumn p_n\MLO(-\la_nt^{\alpha})(\va_n(1) - \psi_n(1)) = 0,
\quad t>0.
$$
Therefore
$$
\sumn \frac{p_n}{\eta^{\gamma}+\la_n}
(\va_n(1) - \psi_n(1)) = 0, \quad \mbox{Re}\, \eta > 0.
$$
Similarly we can obtain
$$
p_m(\va_m(1) - \psi_m(1))  = 0, \quad m\in \N.
$$
Since $p_n\ne 0$ for each $n\in \N$ by {\bf Assumption}, we have
$$
\va_m(1) = \psi_m(1), \quad m\in \N.    \eqno{(3.11)}
$$
Moreover, we assumed
$\va_n(0) = \psi_n(0)=1, \ n\in \N. $

From Lemma 3, we have
$$
\psi_n(x) = \va_n(x) + \int^x_0 K(x,y)\va_n(y) dy, \quad
n\in \N, \, 0<x<1;
$$
and
$$
\frac{d\psi_n}{dx}(x) = \frac{d\va_n}{dx}(x) + K(x,x)\va_n(x)
+ \int^x_0 \ppp_xK(x,y)\va_n(y) dy.
$$
Moreover, we have the boundary condition
$$
\frac{d\psi_n}{dx}(1) + J\psi_n(1)
= \frac{d\va_n}{dx}(1) + H\va_n(1) = 0,
$$
which combining with (3.11) yields
$$
\int^1_0 K(1,y)\va_n(y) dy = 0, \quad n\in \N    \eqno{(3.12)}
$$
and
$$
(J-H+K(1,1))\va_n(1) + \int^1_0 \ppp_xK(1,y)\va_n(y) dy = 0,
\quad n\in \N.                      \eqno{(3.13)}
$$
Since the relation
$$
\lim_{n\to \infty} \int^1_0 \ppp_xK(1,y)\va_n(y) dy = 0
$$
holds, we see
$$
\lim_{n\to \infty} (J-H+K(1,1))\va_n(1) = 0.
$$
Meanwhile, from the asymptotic behavior (2.3) and we note that $\lim_{n\to \infty}\va_n(1) \ne 0$ (e.g., \cite{LS}), we obtain
$$
J-H + K(1,1) = 0.                   \eqno{(3.14)}
$$
Then by (3.13), we readily obtain
$$
\int^1_0 \ppp_xK(1,y)\va_n(y) dy = 0, \quad n\in \N.              \eqno{(3.15)}
$$
In terms of (3.12) and (3.15), since $\{\varphi_n\}_{n\in \N}$ is
an orthogonal basis in $L^2(0,1)$, we see
$$
K(1,y) = \ppp_xK(1,y) = 0, \quad 0<y<1.
$$
Then by the uniqueness of the solution K in Lemma 2, we conclude $K(x,y) =  0$ for $0<y<x<1$.
Therefore (2.5) and (3.14) yield
$$
J=H,
$$
and
$$
j-h + \frac{1}{2}\int^x_0 (q(\xi)-p(\xi))d\xi = 0, \quad 0<x<1.  \eqno{(3.16)}
$$
Hence (3.16) with $x = 0$ yields $j=h$.
Furthermore, differentiate (3.16) and we obtain
$q(x)=p(x)$, $x \in [0,1]$.

Finally, since we already have the uniqueness of $\alpha, p(x), h, H$,
we see that
$$
p_m = q_m, \quad  \quad m\in\N
$$
means
$$
(a, \va_m)=(b, \va_m),  \quad \quad \ m\in\N.
$$
Then we have
$$
a(x)=b(x), \quad \quad x\in [0,1].
$$
Thus the proof of Theorem 1 is complete for the case: $0<\alpha, \beta<1$.

\

{\bf{Case 2:} $\alpha \in (1,2)$.}

Firstly, we have
$$
u(p, h, H, \widetilde{a})(0,t) = u(q, j, J,\widetilde{b})(0,t) \quad 0<t<T,
$$
the analyticity in $t>0$ yields
$$
\sumn p_{n}E_{\alpha,1}(-\la_nt^{\alpha}) + \sumn p^0_{n}tE_{\alpha,2}(-\la_nt^{\alpha})
= \sumn q_{n} E_{\alpha,1}(-\mu_nt^{\alpha}) + \sumn q^0_{n}tE_{\alpha,2}(-\mu_nt^{\alpha}), \quad t > 0.
$$
We take the Laplace
transforms termwise in both side of the above equation to obtain
$$
\sumn \frac{p_n z^{\alpha-1}+ p^0_n z^{\alpha-2}}{z^\alpha +\la_n}
=
\sumn \frac{q_n z^{\alpha-1}+ q^0_n z^{\alpha-2}}{z^\alpha + \mu_n}, \ \ Re \ z>0.    \eqno{(3.17)}
$$
By a similar argument used for $\alpha \in(0,1)$ , we obtain
$$
\la_n = \mu_n, \quad n\geq 1.
$$
Hence we have
\begin{align*}
& \sumn (p_n\MLO(-\la_nt^{\alpha}) + p_n^0 t\MLT(-\la_nt^{\alpha}))
= & \sumn (q_n\MLO(-\la_nt^{\alpha}) + q_n^0 t\MLT(-\la_nt^{\alpha})),
\quad t>0.
\end{align*}
Then taking the Laplace transform, we see
$$
\sumn \frac{p_nz^{\alpha-1} + p_n^0z^{\alpha-2}}{z^{\alpha}+ \la_n}
= \sumn \frac{q_nz^{\alpha-1} + q_n^0z^{\alpha-2}}{z^{\alpha}+ \la_n},
\quad \mbox{Re}\, z > 0.
$$
Similarly to Case 1, let $z \to z_{m, a}^{\pm}$ and set $ z_{m, a}^{\pm}= \la^{\frac{1}{\alpha}}_m e^{\pm \sqrt{-1} \frac{\pi}{\alpha}}$.
Then
$$
p_m(z_{m, a}^{\pm})^{\alpha-1} + p_m^0 (z_{m, a}^{\pm})^{\alpha-2}
= q_m(z_{m, a}^{\pm})^{\alpha-1} + q_m^0 (z_{m, a}^{\pm})^{\alpha-2},
$$
that is,
$$
(p_m-q_m)z_{m, a}^{\pm} + (p_m^0 - q_m^0) = 0, \quad m\in \N.
$$
Hence
$$
p_m = q_m, \quad p_m^0 = q_m^0, \quad m\in \N.    \eqno{(3.18)}
$$

Next, applying the analyticity in $t>0$ to
$u(p, h, H, \widetilde{a})(1,t) = u(q, j, J,\widetilde{b})(1,t)$ $0<t<T$,
we obtain
$$
\sumn p_n\MLO(-\la_nt^{\alpha})(\va_n(1) - \psi_n(1))
+ \sumn p_n^0 t\MLT(-\la_nt^{\alpha})(\va_n(1) - \psi_n(1)) = 0,
\quad t>0.
$$
Therefore
$$
\sumn \frac{p_nz^{\alpha-1}+p_n^0z^{\alpha-2}}{z^{\alpha}+\la_n}
(\va_n(1) - \psi_n(1)) = 0, \quad \mbox{Re}\, z > 0.
$$
Similarly we can obtain
$$
(\va_m(1) - \psi_m(1))(p_m z_{m, a}^{\pm} + p_m^0) = 0, \quad m\in \N.
$$
By {\bf Assumption}, we have $p_m\ne 0$ or $p_m^0 \ne 0$.  Hence,
$$
p_m z_{m,a}^+ + p_m^0 \ne 0 \quad \mbox{or} \quad p_m z_{m, a}^- + p_m^0 \ne 0.
$$
Therefore,
$$
\va_n(1) = \psi_n(1), \quad n\in \N.    \eqno{(3.19)}
$$
By Lemma 3 and $\va_n(0) = \psi_n(0) = $ for $n\in \N$, we have
$$
\psi_n(x) = \va_n(x) + \int^x_0 K(x,y)\va_n(y) dy, \quad
n\in \N, \, 0<x<1;
$$
and
$$
\frac{d\psi_n}{dx}(x) = \frac{d\va_n}{dx}(x) + K(x,x)\va_n(x)
+ \int^x_0 \ppp_xK(x,y)\va_n(y) dy.
$$
Then similarly to the case $\alpha \in(0,1)$, we can obtain
$$
J=H, \quad j=h
$$
and
$$
q(x)=p(x), \ \ x \in [0,1].
$$

Finally, since we already have the uniqueness of $\alpha, p(x), h, H$,
then (3.18) imply that
$$
(a, \va_m)=(b, \va_m) \quad and \quad (a^0, \va_m)=(b^0, \va_m), \ m\in\N.
$$
Hence we have
$$
a(x)=b(x),\quad and \quad a^0(x)=b^0(x).
$$
Thus the proof of Theorem 1 is complete.

\section{Proof of Theorem 2}
\sloppy{}

In this part, we give the proof of Theorem 2 for problem (1.6)-(1.9).
 Based on the uniqueness result in Theorem 1
and using the Duhamel principle and the Titchmarsh convolution theorem,
we reduce the proof of Theorem 2 to Theorem 1 as follows.
\\
\\
{\bf Proof of Theorem 2.} We consider Theorem 2 separately in following three cases.

{\bf Case 1: $0<\alpha<1$.}
From the fractional Duhamel principle in e.g., \cite{LRY,LYZZ}, we have
$$
\frac{1}{\Gamma(1-\alpha)}\int^t_0 (t-s)^{-\alpha}
\widetilde{u}(\ell,s) ds
= \int^t_0 \theta(t-s) u (\ell, s) ds \eqno{(4.1)}
$$
for $\ell=0,1$ and $0<t<T$.\\
Therefore,
$\widetilde{u}(p, h, H)(\ell, t) = \widetilde{u}( q, j, J)(\ell, t)$, $0<t<T$, $\ell=0,1$ yield
$$
\int^t_0 \theta(t-s) \left(u( p, h, H)(\ell, t) - u( q, j, J)\right)(\ell,s) ds = 0, \quad \ell=0,1, \,
0<t<T. \eqno{(4.2)}
$$
The Titchmarsh convolution theorem (\cite{ECT}) implies the existence
of $T_1, T_2 \geq 0$ satisfying $T_1 + T_2 \geq T$ such that $\theta(t) = 0$ for
$t\in(0, T_1)$ and $ u( p, h, H)(\ell, t) - u(q, j, J)(\ell,t)= 0$ for $t \in [0, T_2]$.
However, since $\theta\in C^1(0,T), \not\equiv 0$, this implies that $\theta(t) > 0$
a.e. in $(0, T)$. As a result, we obtain $T_1 = 0 $
and thus $ T_2 = T$, that is,
$$
u( p, h, H)(\ell, t) = u( q, j, J)(\ell,t), \quad \ell=0,1, \, 0<t<T. \eqno{(4.3)}
$$
Thus the proof is reduced to Theorem 1, and we omit further details.

\

{\bf Case 2: $\alpha=1$.}
The Duhamel principle in \cite{LYZZ} yields
$$
\widetilde{u}(p, h, H)(\ell, t) = \int^t_0 \theta(t-s) u(p, h, H)(\ell, s) ds.
$$
for $\ell=0,1$ and $0<t<T$.
Therefore, using the same arguments in the proof of case 1, and the result of Theorem 1,
we can argue to complete the proof.

\

{\bf Case 3: $1<\alpha<2$.}
From the fractional Duhamel principle (\cite{GYM}), we have
$$
\frac{1}{\Gamma(2-\alpha)}\int^t_0 (t-s)^{1-\alpha}\widetilde{u}(p, h, H)(\ell, s) ds
= \int^t_0 \theta(t-s) u(p, h, H)(\ell, s) ds
$$
for $\ell=0,1$ and $0<t<T$.
Therefore, similarly to the proof in case 1,
we can complete the proof.
Thus the proof of Theorem 2 is complete.

\section*{Appendix. Proof of Lemma 4}
\sloppy{}

The proof is similar to Step I of the proof of Theorem 1 in Section 3.
Indeed we consider only the case $1<\alpha<2$, because the case of
$0<\alpha \le 1$ is similar.
Then by the same way as Case 3 in Step I in Section 3 for obtaining (3.7),
we have
$$
\sumn \frac{p_nz+p_n^0}{z^{\alpha}+\la_n} = 0, \quad z \in Q_0
:= \{ z\in \C;\, z\ne 0, \, \vert arg \ z \vert < \frac{\pi}{2}\}
\cap (\C \setminus \{ z_{m,\alpha}^{\pm}\}_{m\in \N}),
$$
where we set $z_{m,\alpha}^{\pm} := \la_m^{\frac{1}{\alpha}}e^{\pm \sqrt{-1}
\frac{\pi}{\alpha}}$.
For arbitrary $n_0 \in \N$, setting $J_0(z):= - \sum_{n\ne n_0}
\frac{p_nz+p_n^0}{z^{\alpha}+\la_n}$, we see
$$
p_{n_0}z + p_{n_0}^0 = (z^{\alpha}+\la_{n_0})J_0(z), \quad z\in Q_0.
$$
We set $z_{n_0} := z_{n_0,\alpha}^+ = \la_{n_0}^{\frac{1}{\alpha}}
e^{\sqrt{-1}\frac{\pi}{\alpha}}$.  Then $z_{n_0}^{\alpha} = -\la_{n_0}$ and
$\lim_{z\in Q_0, z\to z_{n_0}} J_0(z) = J_0(z_{n_0})$ by noting
$1<\alpha<2$.  Therefore $p_{n_0}z_{n_0} + p_{n_0}^0 = 0$.
Since $p_{n_0}, p_{n_0}^0 \in \R$ and $z_{n_0} \not\in \R$ by $1<\alpha<2$,
we obtain $p_{n_0} = p_{n_0}^0 = 0$.  Since $n_0\in \N$ is arbitrary,
we see that $p_n = p_n^0 = 0$ for all $n\in \N$.  By (3.2) we reach
$\UUU=0$ in $(0,1) \times (0,T)$.  Thus the proof of Lemma 4 is complete.

\section*{Acknowledgements}
\sloppy{}

The first author thanks the China Scholarship Council for their support.
The second author was supported by Grant-in-Aid for Scientific Research (S)
15H05740 and Grant-in-Aid (A) 20H00117 of
Japan Society for the Promotion of Science and
by The National Natural Science Foundation of China
(no. 11771270, 91730303).
This paper has been supported by the RUDN University
Strategic Academic Leadership Program.

\end{document}